\documentstyle[amssymb,amsfonts,amscd,11pt]{amsart}

\newtheorem{theorem}{Theorem}[section]
\newtheorem{proposition}[theorem]{Proposition}

\newtheorem{corollary}[theorem]{Corollary}

{\theoremstyle{break}

}

\newtheorem{definition}[theorem]{Definition}

\newcommand{\vrpf}[2]{\begin{picture}(#1,1)
\put(0,1){\vector(1,0){#1}} \put(0,2.5){\makebox(#1,0)[b]{$#2$}}
\end{picture}}
\newcommand{\vupf}[2]{ \Vupf{#1} \makebox(0,#1)[l]{$\,#2$}}
\newcommand{\Vupf}[1]{\begin{picture}(0.5,#1)
                 \put(0,#1){\vector(0,-1){#1}}\end{picture}}
\newcommand{\vupflr}[3]{\makebox(0,#1)[r]{$#2\;$}
\vupf{#1}{#3}}

\newcommand{\cL}{{\mathcal L}}
\newcommand{\cO}{{\mathcal O}}
\newcommand{\cX}{{\mathcal X}}
\newcommand{\Gal}{{\mbox{\rm Gal}}}
\newcommand{\ms}{\medskip}

\newcommand{\ov}{\overline}
\newcommand{\bbbc}{{\mathbb C}}
\newcommand{\bbbf}{{\mathbb F}}
\newcommand{\bbbp}{{\mathbb P}}
\newcommand{\bbbq}{{\mathbb Q}}
\newcommand{\bbbr}{{\mathbb R}}
\newcommand{\bbbg}{{\mathbb G}}
\newcommand{\bbbz}{{\mathbb Z}}
\newcommand{\ba}{\begin{array}}
\newcommand{\ea}{\end{array}}
\newcommand{\univ}{{\mbox{\rm\scriptsize univ}}}
\newcommand{\id}{{\mbox{\rm id}}}
\newcommand{\Aut}{{\mbox{\rm Aut}}}
\newcommand{\Gr}{{\mbox{\rm Gr}}}
\newcommand{\Sm}{{\mbox{\rm Sm-et}}}
\newcommand{\tr}{{\mbox{\rm tr}}}
\newcommand{\rk}{{\mbox{\rm rk}}}
\newcommand{\SL}{{\mbox{\rm SL}}}
\newcommand{\pr}{{\mbox{\rm pr}}}

\newcommand{\beqn}{\begin{eqnarray*}}
\newcommand{\eeqn}{\end{eqnarray*}}
\begin{document}

\title[Moduli Stacks of Vector Bundles and Frobenius Morphisms]
{Moduli Stacks of Vector Bundles and Frobenius morphisms}

\author[Frank Neumann]
{Frank Neumann}
\address{Department of Mathematics and Computer Science\\
University of Leicester\\
University Road, Leicester LE1 7RH, England, UK}
\email{fn8@@mcs.le.ac.uk}

\author[Ulrich Stuhler]
{Ulrich Stuhler}
\address{Mathematisches Institut\\
Georg-August-Universit\"at G\"ottingen\\
Bunsenstr. 3-5, D-37073 G\"ottingen, Germany}
\email{stuhler@@uni-math.gwdg.de}




\begin{abstract}
{We describe the action of the different Frobenius morphisms on
the cohomology ring of the moduli stack of algebraic vector
bundles of fixed rank and determinant on an algebraic curve over a
finite field in characteristic $p$ and analyse special situations
like vector bundles on the projective line and relations with
infinite Grassmannians.}
\end{abstract}

\maketitle

\vspace*{-0.6cm}

\section*{Introduction}


\ms Suppose, $F:X\to X$ denotes the geometric Frobenius
endomorphism of a smooth, projective algebraic curve of genus $g$
over the finite field $\bbbf_q$ of $q=p^s$ elements of
characteristic $p>0$. It is well known and has been studied quite
often, that the pullback operation on vector bundles on $X$,
induced by $F$, does not necessarily respect the stability of the
bundle. If the genus of the curve $X$ is greater than one, not
much is known about this situation.

\ms In particular, the pullback operation $F^{\ast}$ does not in
general induce a morphism of the coarse moduli scheme of bundles
in the sense of Narasimhan and Seshadri, but only a rational map.
It is only recently, that Y. Laszlo and C.\ Pauly have been able
to write down such a rational map explicitly in the special case
of rank two vector bundles for $p=2,\,  g=2$, making use of the
explicit knowledge of the coarse moduli scheme of (semi) stable
bundles of rank $2$ and degree $0$ in the case $g=2$
(\cite{Las-Pau1} and  \cite{Las-Pau2} and for more results
\cite{Jo-Ra}).

\ms It is therefore natural to consider the action of $F^{\ast}$
or $\overline{F}^{\ast}$, extending $F$ to the algebraic closure
$\overline{\bbbf}_q$ of $\bbbf_q$, on the stack of bundles of rank
$n$ with trivial determinant as it is done in this paper.

\ms Our contribution here consists of the observation (Prop. 2.4)
that despite the fact that we do not know much about the action of
$F^{\ast}$, it is possible to evaluate the action of $F^{\ast}$ or
$\overline{F}^{\ast}$ on the cohomology ring of the moduli stack
of vector bundles.

\ms Of course, the next question is about the existence of a
Lefschetz trace formula in this context, and we show by some easy
examples that this can not work, at least not in a naive sense. In
fact, already the case of the projective line $X=\bbbp^1$ with
$g=0$ and for simplicity, $n=2$, shows that there have to be
modifications for such a trace formula taking care of fixed points
at infinity of the stack of bundles. Formally, one can force
convergence of the corresponding infinite sums (the stack of
bundles has cohomology in infinitely many dimensions!) by twisting
with a high enough power of the arithmetic Frobenius of the stack
of bundles. Nevertheless this convergent expression is not
directly related to the set of fixed points, as is shown already
by the example above. It should be mentioned at this point that a
Lefschetz trace formula for stacks works quite well for actions of
the powers of the arithmetic Frobenius of the stack of bundles
alone (without the pullback $\overline{F}^{\ast}$) as is discussed
in \cite{Be}, \cite{Lau-Mo-Bai}.

\ms The organisation of this paper is as follows: the first
section contains some general material concerning stacks of vector
bundles and describes the different Frobenius morphisms occuring.

\ms Section two contains the cohomological computations and
evaluates in particular the action of $\overline{F}^{\ast}$ on the
cohomology of the stack of bundles. Concerning the cohomology of
the stack of bundles we have made use of work of Behrend \cite{Be}
and in particular of the nice diploma thesis of J. Heinloth
\cite{Hei}. The procedure followed in their papers is parallel to
the fundamental work of Atiyah and Bott \cite{At-Bo}, but whereas
Atiyah and Bott work in a differential geometrical context the
procedure here is in the language of algebraic geometry and
stacks. The decisive point is the cohomology of the gauge group
which has to be treated here in a different way using older
results of Harder and Narasimhan on the Tamagawa number of the
special linear group. For the convenience of the reader as part of
the work \cite{Hei} is unpublished, we will repeat here some of
the arguments of \cite{Hei} with hints to the literature and in
particular to \cite{Bi-Ghi-Lat}, where similar results are shown
in a slightly different way.

\ms The third section gives in part one the discussion of the
example $X=\bbbp^1$ of the projective line and bundles of rank
$n=2$. The interesting situation is here, that thanks to
Grothendieck's splitting theorem for vector bundles, the operation
of $\overline{F}^{\ast}$ on the stack is completely explicit. As
is mentioned already, a trace formula in a naive sense does not
hold. On the other hand we have the definite impression, that, if
one could make sense of our approach in this case, this should
generalize to the cases of arbitrary genus $g$, where not much is
known about $\overline{F}^{\ast}$.

\ms Part two of section three contains some further computations
of the cohomology of moduli stacks of bundles and discusses the
description of the stack of bundles as some kind of double
quotient using the uniformisation theorem of Drinfeld and Simpson
\cite{Dri-Sim} (see also \cite{Tel}). The ingredients of this
description are the (infinite) affine Grassmannians from the
theory of loop groups and an algebraic version of the gauge group.
We indicate here a computation of the cohomology of both objects,
evaluating along this a certain Leray spectral sequence. These
considerations and the double coset description might be useful
for extending our approach to a compactified situation.

\ms

{\it Acknowledgements.} The first author likes to thank all the
people at the Centre de Recerca Matem\`atica, Bellaterra, Spain
and the members of the Barcelona Algebraic Topology Group for much
support and a wonderful stay, while holding a post-doc research
grant as part of the EU research training network Modern Homotopy
Theory. He also likes to thank the people at the Mathematisches
Institut of the Georg-August-Universit\"at G\"ottingen, where this
work was initiated. The second author has reported on these
questions in February 2001, at the Tata Institute, Mumbai, India
and would like to thank the people there and in particular N.
Nitsure and S. Ramanan for this beautiful stay and the possibility
to speak about these things. Both authors also like to thank the
referee for many useful suggestions and remarks.

\section{Vector bundles on curves and Frobenius morphisms}

\ms Let $X$ denote a smooth, complete and irreducible algebraic
curve of genus $g$ over the finite field $\bbbf_q$ with $q$
elements. For any scheme $S$ over $\bbbf_q$, $\tilde{\cL} (S)$
denotes the following category: the objects of $\tilde{\cL} (S)$
are the rank $n$ vector bundles over the scheme $X\times_{\bbbf_q}
S$, the morphisms of $\tilde{\cL} (S)$ are just isomorphisms
between such vector bundles.

\ms Additionally, we consider the category $\cL (S)$, whose
objects are pairs $(E,\delta)$, where $E$ ist a rank $n$ vector
bundle over $X\times_{\bbbf_q} S$ and $\delta$ is an isomorphism
of line bundles
$$\delta : \det (E) \stackrel{\sim}{\to} \cO_{X\times_{\bbbf_q}S}$$
of the determinant bundle of $E$ with the trivial line bundle. We
always identify vector bundles and their locally free module
sheaves of sections.
 The morphisms $u:(E,\delta) \to (E', \delta')$ in $\cL(S)$ are  isomorphisms
of vector bundles
$$
u:E \stackrel{\sim}{\longrightarrow} E',
$$
inducing commutative diagrams
$$
\begin{array}{ccc}
\det(E)&\vrpf{10}{\det(u)}&\det(E')\\
\vupflr{10}{\delta}{\sim}& &\vupflr{10}{\sim}{\delta'}\\
\cO_{X\times_{\bbbf_q}S}&\vrpf{10}{\id}&\cO_{X\times_{\bbbf_q}S}
\end{array}
$$

It is well known, that the functor
$$
S/ \bbbf_q \mapsto \tilde{\cL} (S)
$$
on the category of schemes over $\bbbf_q$, defines an algebraic
stack in the sense of M. Artin (\cite{Lau-Mo-Bai}, th\'eor\`eme
4.6.2.1). Due to the fact that $X$ is a curve over $\bbbf_q$, it
follows, that $\cL$ is a smooth algebraic stack. (\cite{Fal1},
\cite{Hei}, 2.1.3, p.39).

\ms Besides $\tilde{\cL}$, we have also the stack $\cL$, given by
the functor $S/\bbbf_q \mapsto \cL (S)$, which is again smooth.

\ms We will consider also the corresponding functors and stacks
over $\ov{\bbbf}_q$, denoted by $\ov{\tilde{\cL}}$ and $\ov{\cL}$
respectively. From the definitions it follows immediately that
there is a canonical isomorphisms of stacks
\begin{eqnarray*}
\tilde{\cL} \times_{\bbbf_q} \ov{\bbbf}_q & \widetilde{\to} &
\ov{\tilde{\cL}},\\
\cL \times_{\bbbf_q}\ov{\bbbf}_q & \widetilde{\to} & \ov{\cL}.
\end{eqnarray*}

\ms On $X\times_{\bbbf_q} \tilde{\cL}$ we have the universal
vector bundle $E^{\univ}$, such that for any scheme $S$ over
$\bbbf_q$ a morphism $\varphi:S\to \tilde{\cL}$ induces the vector
bundle $(\id_X \times \varphi)^{\ast} (E^{\univ})$ on
$X\times_{\bbbf_q} S$, which is exactly the element of
$\tilde{\cL} (S)$ given by $\varphi$.

\ms In the same way, on $X\times_{\bbbf_q} \cL$, we have a
universal pair $(E^{\univ}, \delta_{\univ})$, where
$\delta_{\univ}:\det (E^{\univ}) \stackrel{\sim}{\to}
\cO_{X\times_{\bbbf_q} \cL}$. Similar considerations are true over
the algebraic closure  $\ov{\bbbf}_q$.

\ms Following \cite{At-Bo}, \cite{Fal1}, \cite{Hei}, we will
consider various open substacks of $\tilde{\cL}$ resp. $\cL$ and
similarly for $\ov{\tilde{\cL}}$ resp. $\ov{\cL}$.

\ms For an arbitrary vector bundle $E$ over an algebraic curve $X$
(over a field), $\mu (E)$ denotes the quotient
$$
\mu (E) := \frac{\deg (E)}{\rk(E)},
$$
where $\deg (E)$ is the degree of the bundle $E$, $\rk(E)$ its
rank.

\ms {\definition A vector bundle $E$ is stable (resp. semistable)
if for all proper subsheaves (resp. all subsheaves) $E'$ of $E$
the inequality $\mu (E') < \mu (E)$ (resp. $\mu (E') \le \mu (E)$)
holds.}

\bigskip
As is well known, any vector bundle $E$ on a curve $X$ has a
canonical filtration, the {\it Harder-Narasimhan filtration}:
$$
E_0 = (0) \subset E_1 \subset \ldots \subset E_r =E,
$$
such that the following properties hold:
\begin{enumerate}
\item[i)] The successive quotients $E_i/ E_{i-1}$ are all semistable.
\item[ii)] Given $E_i, E_{i+1}$ is a maximal subbundle of $E$ with the
additional property that
$$
\mu (E_{i+1}/E_i) \ge \mu (E')
$$
holds for all subbundles $E' \subset E/E_i$. Given $E_i$, this
determines $E_{i+1}$ uniquely.
\item[iii)] Dually, given $E_i$, then $E_{i-1}$ is the smallest subbundle of
$E_i$ with the property, that
$$
\mu (E_i/E_{i-1}) \le \mu (E'')
$$
for all quotient bundles $E_i/E''$.
\end{enumerate}

The Harder-Narasimhan filtration is uniquely determined by the
above properties. We associate now with this filtration the
polygonal function with $\rk(E)=n$
$$
p_E : [0, \rk(E)] \to \bbbr,
$$
given as follows: $p_E (\rk(E_i))=\deg (E_i)$ for $i=0, \ldots,
r$. At all other values $x \in [0, \rk(E)]$, the value $p_E(x)$ is
obtained by linear interpolation using the values defined before.
Obviously $p_E$ is a concave function, that is, slopes are
decreasing with growing $x$.

\ms For any piecewise linear function $p:[0, \rk (E)] \to \bbbr$,
satisfying $p(\rk(E)) = \deg (E) =0$ (we assume, that $\det (E)
\cong \cO_X$), we define substacks $\cL_{\le p}$ resp. $\cL_{<p}$.

\ms {\definition For any scheme $S$ over $\bbbf_q$, $\cL_{\le p}
(S)$ is the subcategory of $\cL(S)$ of pairs $(E,\delta)$, $E$ a
bundle over $X\times_{\bbbf_q} S$ of rank $n$, $\delta : \det (E)
\stackrel{\sim}{\to} \cO_{X\times_{\bbbf_q}S}$ a trivialisation of
the determinant bundle, such that for all closed points $s \in S$,
the Harder-Narasimhan polygon $p_{E_s}$ of $E\times_{\bbbf_q}s$
satisfies $p_{E_s} \le p$. Similarly, $\cL_{< p} (S)$ denotes the
subcategory of pairs $(E, \delta)$, such that $p_{E_s} < p$ holds
for all closed points $s \in S$.}

\ms {\proposition The stacks $\cL_{<p} \hookrightarrow \cL$ resp.
$\cL_{\le p} \hookrightarrow \cL$ are open substacks of $\cL$.
Similarly $\cL_{<p} \hookrightarrow \cL_{\le p}$ is an open
substack. $\cL_p$ denotes the reduced closed substack $(\cL_{\le
p} \setminus L_{<p})$.}

\ms {\bf Proof.} {The proposition follows by using the
semicontinuity theorem. For the details see for example
\cite{Hei}, 2.1.10.} \qed

\ms {\bf Remark.} Similarly we have the open substacks
$\ov{\cL}_{\le p} = \cL_{\le p} \times_{\bbbf_q} \ov{\bbbf}_q$ and
$\ov{\cL}_{< p} = \cL_{<p} \times_{\bbbf_q} \ov{\bbbf}_q$ of
$\ov{\cL}$ as well as the closed substack $\ov{\cL}_p= \cL_{p}
\times_{\bbbf_q} \ov{\bbbf}_q$ of $\ov{\cL}$.

\bigskip
We recall for completeness

\ms {\definition $(E^{\univ}, \delta^{\univ})$ denotes the
universal rank $n$ vector bundle over $X\times_{\bbbf_q} \cL$ with
trivial determinant $\delta^{\univ} :\det (E) \stackrel{\sim}{\to}
\cO_{X\times_{\bbbf_q} \cL}$.}

\bigskip
We consider the following Frobenius morphisms in this context.
First we have the geometric Frobenius endomorphism
$$
F_X :(X, \cO_X) \to (X, \cO_X),
$$
$$
F_X = \id\, \,  \mbox{\rm (on X)}\, \, , F^{\ast}_X (f) = f^q \,
\, \mbox{\rm (on sections of $\cO_X$).}
$$

\ms Let $\ov{F}:= F_X \times_{\bbbf_q} \mbox{\rm
id}_{\ov{\bbbf}_q}$ denote the extension of $F_X$ to an
endomorphism of $X\times_{\bbbf_q} \ov{\bbbf}_q$ over
$\ov{\bbbf}_q$.

\ms If $S/\ov{\bbbf}_q$ denotes a scheme over $\ov{\bbbf}_q$, one
has the pullback operation
\begin{eqnarray*}
\ov{\cL} (S) & \to & \ov{\cL} (S)\\
(E, \delta) & \mapsto & \big(\ov{F}^{\ast} (E), \ov{F}^{\ast}
(\delta)\big).
\end{eqnarray*}
This induces an endomorphism of the stack $\ov{\cL}$,
$$
\varphi : \ov{\cL} \to \ov{\cL}.
$$
{\bf Remark.} As already outlined in the introduction, the nature
of $\varphi$ is very mysterious. In particular it does not respect
the open substacks $\ov{\cL}_{\le p} \hookrightarrow \ov{\cL}$ or
$\ov{\cL}_{<p} \hookrightarrow \ov{\cL}$.

\ms Besides $\varphi$ we have also the geometric Frobenius
morphism
$$F_{\cL} : ( \cL, \cO_{\cL}) \to (\cL, \cO_{\cL})$$

and its extension $F_{\cL} \times_{\bbbf_q} \mbox{\rm
id}_{\ov{\bbbf}_q}$ as an endomorphism of $\ov{\cL}$ over
$\ov{\bbbf}_q$. Of course this endomorphism respects the open
substacks $\cL_{\le p} \hookrightarrow \cL$, $\cL_{<p}
\hookrightarrow \cL$ and their extensions over $\ov{\bbbf}_q$.\\
Actually we prefer to work in the next chapters with the
arithmetic Frobenius morphism on the smooth-\'etale site $(\cL
\times_{\bbbf_q}\ov{\bbbf}_q)_{\tiny{\Sm}}$, given as $\psi:=
\id_{\cL} \times \, \mbox{\rm Frob} \big(\ov{\bbbf}_q /
\bbbf_q\big)$ which acts as an inverse to the action of $F_{\cL}
\times_{\bbbf_q} \id_{\ov{\bbbf}_q}$ on the smooth-\'etale site
$(\cL \times_{\bbbf_q} \ov{\bbbf}_q)_{\tiny{\Sm}}$.

\ms {\proposition There is a canonical isomorphism
$$\big(\ov{F}_X \times \id_{\ov{\cL}}\big)^{\ast}\big(E^{\univ}\big)
\cong \big(\id_{\ov{X}} \times
\varphi\big)^{\ast}\big(E^{\univ}\big).$$}

{\bf Proof.} This is immediate. For a stack $T/\ov{\bbbf}_q$ a
rank $n$ vector bundle $E$ on $\ov{X} \times T$ together with an
isomorphism
$$
\det(E) \stackrel{\sim}{\to} \cO_{\ov{X}\times_{\bbbf_q} T}
$$
is given by a morphism $u:T \to \ov{\cL}$, such that $E\cong
(\id_{\ov{X}}\times u)^{\ast}(E^{\univ})$ and similarly for the determinant.\\
\noindent Applying this to the vector bundle $\big(\ov{F}_X
\times_{\ov{\bbbf}_q} \id_{\ov{\cL}}\big)^{\ast} (E^{\univ})$,
which defines $\varphi : \ov{\cL} \to \ov{\cL}$, and the
proposition follows. \qed

\section{Cohomology of the moduli stacks of vector bundles}

\ms In this section we describe the cohomology of the moduli stack
of vector bundles. We collect the results from various different
treatments in the literature. (see \cite{Hei}, \cite{Be},
\cite{Bi-Ghi-Lat} and \cite{At-Bo}).

\ms Let $X$ be an algebraic curve over the field $\bbbf_q$. We
have the following well known description of the $l$-adic
cohomology ring of the curve $\ov{X} = X\times_{\bbbf_q}
\ov{\bbbf}_q$ over the algebraic closure $\ov{\bbbf}_q$, where
$\ov{\bbbq}_{l}$ denotes the algebraic closure of $\bbbq_{l}$,
namely \beqn
H^0 (\ov{X}; \ov{\bbbq}_l) & = & \ov{\bbbq}_l \cdot 1\\
H^1 (\ov{X}; \ov{\bbbq}_l) & = & \mathop{\oplus}\limits^{2g}_{i=1}
\ov{\bbbq}_l
\cdot \alpha_i\\
H^2 (\ov{X}; \ov{\bbbq}_l) & = & \ov{\bbbq}_l [\ov{X}], \eeqn
where $[\ov{X}]$ denotes the orientation class of
$\ov{X}/\ov{\bbbf}_q$ and where we can assume, that the
$\{\alpha_i : i=1, \ldots, 2g\}$ are eigenclasses under the action
of the Frobenius morphism. The Frobenius morphism acts semisimple
on the etale cohomology of a curve (see \cite{Mu}, p.203).

\ms More precisely, for $\ov{F}=F_X \times_{\bbbf_q}
\id_{\ov{\bbbf}_q}$ as before, we have
 explicitely that
\beqn
\ov{F}^{\ast} (1) & = & 1\\
\ov{F}^{\ast} (\alpha_i) & = & \lambda_i \alpha_i \quad (i=1,
\ldots, 2g), \eeqn where $\lambda_i \in \ov{\bbbq}_l$ is
algebraic, such that $|\lambda_i| = q^{\frac12}$ for any possible
embedding of $\lambda_i$ into the complex numbers $\bbbc$. Finally
$\ov{F}^{\ast} ([\ov{X}]) = q \cdot [\ov{X}]$ for the orientation
class $[\ov{X}]$.

\ms The definition of the $l$-adic cohomology of an algebraic
stack $\cX$ can be found in \cite{Lau-Mo-Bai}, chapter 12. In
\cite{Be}, \cite{Hei} the authors work with the the smooth
topology on algebraic stacks. We will follow here the treatment of
\cite{Lau-Mo-Bai}. For an algebraic stack $\cX$ over a scheme $S$,
$\cX/S$, we consider the smooth-\'etale site $\Sm (\cX/S)$
(\cite{Lau-Mo-Bai}, Def. 12.1.) The underlying category of $\Sm
(\cX/S)$ has as objects the smooth $1$-morphisms of $S$-algebraic
stacks $u:U\to\cX$, where $U$ is an algebraic space over $S$. A
morphism is a pair $(\varphi,\alpha):(U,u)\to (V,v)$, where
$\varphi:U\to V$ is a morphism of $S$-algebraic spaces and
$\alpha$ is a 2-isomorphism $\alpha:u\tilde{\to}(u\circ\varphi)$.
The topology on this category is generated by the pretopology of
families of morphisms
$$
((\varphi_i,\alpha_i):(U_i,u_i)\to(U,u))_{i\in I},
$$
where the 1-morphism $\coprod\limits_{i\in
I}\varphi:\coprod\limits_{i\in I}
U_i\to U$ is \'etale and surjective.\\
The $l$-adic cohomology of a stack $\cX$ is defined with respect
to this site $\cX_{\tiny{\Sm}}$ as
$$
H^\bullet(\cX;\bbbq_l)=(\lim\limits_{\overleftarrow{n}}H^\bullet(\cX;\bbbz/l^n\bbbz))\otimes_{\bbbz_l}\bbbq_l.
$$

Basically, the method is to take a covering  $X \to \cX$ (over
$\ov{\bbbf}_q$) of the stack $\cX$ by a scheme $X$ and to consider
the cohomology of the associated simplicial scheme

$$\renewcommand{\arraystretch}{-1.7}
[X/\cX]:=[X\begin{array}{c}\rightarrow\\ \gets\\ \gets\end{array}
X\times_{\cX} X
\begin{array}{c}\gets\\
\gets\\
\gets\\
\rightarrow\\
\rightarrow\end{array} X\times_{\cX} X \times_{\cX} X \ldots]
$$
and an appropriate smooth-\'etale site over it. For details, see
chapters 12 and 18 of the book \cite{Lau-Mo-Bai} or \cite{Tel}.

\ms We describe now in detail the steps of the computation of the
cohomology $H^\bullet(\ov{\cL};\bbbq_l)$ using the treatment of
\cite{Bi-Ghi-Lat} and comparing it with \cite{Be} and \cite{Hei}.
The definition of $l$-adic cohomology above causes no problems for
the algebraic stack of bundles $\ov{\cL}$, as $\ov{\cL}$ is the
inductive limit of the open substacks $\ov{\cL}_{\le p}$ and the
cohomology of this inductive system over the Harder-Narasimhan
polygons $p$ is constant for large $p$ as will be discussed
below.\\

\ms {\bf 1)} As described in section 1, Definition 1.3. and
Proposition 1.4., we have upon fixing a Harder-Narasimhan polygon
$p$, the open substacks, $\ov{\cL}_{<p}\hookrightarrow \ov{\cL}$
as well as $\ov{\cL}_{\le p}\hookrightarrow \ov{\cL}$. Similarly
$\ov{\cL}_{<p}\hookrightarrow \ov{\cL}_{\le p}$ is an open
substack, $\ov{\cL}_p$ denotes the reduced closed substack
$(\ov{\cL}_{\le p}\setminus \ov{\cL}_{<p})$.\\
In particular, because $X$ is one-dimensional, $\ov{\cL}_{\le p},
\ov{\cL}_{<p},\ov{\cL}_p$ are smooth algebraic stacks and
$\ov{\cL}_p\hookrightarrow \ov{\cL}_{\le p}$ is a smooth pair.
This can be found in a similar way also in \cite{Bi-Ghi-Lat},
section 5.

\ms {\bf 2)} For the smooth pair $i:\ov{\cL}_p\hookrightarrow
\ov{\cL}_{\le p}$ one has the tangent stacks $T(\ov{\cL}_p)$ and
$T(\ov{\cL}_{\le p})$, (\cite{Lau-Mo-Bai}, 17.11 -- 17.15 and
17.16 -- 17.17), as well as the normal bundle
$$
N(\ov{\cL}_{\le p}/\ov{\cL}_p)=[i^*T(\ov{\cL}_{\le
p})/T(\ov{\cL}_p)].
$$
The normal bundle has the following explicit description in terms
of the
universal bundle $(E^{\rm univ},\delta^{\rm univ})$ (see Definition 1.5).\\
Upon restricting $(E^{\univ},\delta^{\univ})|_{\ov{\cL}_p}$, one
has the canonical sequence of endomorphism bundles
$$
0\to \mbox{\rm End}^{(0)}_{\rm
filt}(E^{\univ}|\ov{\cL}_p)\to\mbox{\rm
End}^{(0)}(E^{\univ}|\ov{\cL}_p)\to \widetilde{\mbox{\rm
End}}^{(0)}(E^{\univ}|\ov{\cL}_p)\to 0,
$$
where ``$(0)$'' denotes the ``trace=0''endomorphism bundles and
``filt'' specifies the endomorphisms respecting the
Harder-Narasimhan filtration (fibrewise). It is easy to see, that
there is an isomorphism of bundles
$$
N(\ov{\cL}_{\le p}/\ov{\cL}_p)=R^1\pr_*(\widetilde{\mbox{\rm
End}}^{(0)}(E^{\univ}|\ov{\cL}_p))
$$
(see \cite{Hei}, 1.3, p. 37 and \cite{Bi-Ghi-Lat}, Prop. 5.2, in
particular 5.2.(3)).

\ms {\bf 3)} We have a morphism of stacks
$$
f:\ov{\cL}_p\to
(\prod\limits^l_{i=1}\ov{\tilde{\cL}}(X;n_i;d_i))^{(0)},
$$
which associates with a vector bundle of a fixed Harder-Narasimhan
filtration of type $p$ the semistable subquotients occuring in the
filtration. ``$(0)$'' here denotes the bundles equipped with a
trivialisation of the total determinant. This morphism of
algebraic stacks is not representable. Nevertheless it is not
difficult to conclude, using the associated simplicial schemes,
that $f$ is cohomologically acyclic inducing in particular an
isomorphism on cohomology (with constant coefficients), (see
\cite{Hei}, 2.1.1. and \cite{Bi-Ghi-Lat}, Prop. 7.1, 7.2)).

\ms {\bf 4)} For the smooth pair $(\ov{\cL}_{\le p},\ov{\cL}_p)$
one has a Gysin sequence
\begin{eqnarray*}
\ldots
&\to&H^{\bullet-2c}(\ov{\cL}_p;\bbbz/l^n\bbbz(c))\stackrel{i_0}{\longrightarrow}H^\bullet(\ov{\cL}_{\le
p};\bbbz/l^n\bbbz)\\
 &\to& H^\bullet(\ov{\cL}_{<p};\bbbz/l^n\bbbz)\longrightarrow H^{\bullet
+1-2c}(\ov{\cL}_p;\bbbz/l^n\bbbz(c)) \to \ldots\quad
 \end{eqnarray*}
where $c$ is the codimension of the stack $\ov{\cL}_p$ in the
stack $\ov{\cL}_{\le p}$.

\ms {\bf Proof.} This is \cite{Be}, Prop. 2.1.2 and Corollary
2.1.3.

\ms {\bf 5)} The composition $i^*i_*$ is given as the cup-product
$$
i^*i_*(x)=c_r(N^*)\cup x,
$$
where $N=N(\ov{\cL}_{\le p}/\ov{\cL}_p)$ is the normal bundle
computed above, $N^*$ the dual bundle, $c_r(N^*)$ is the top Chern
class of the bundle $N^*$.

\ms {\bf Proof.} For schemes this is \cite{Sem}, Expos\'e VII,
Th\'eor\`eme 4.1., p. 299. The case of stacks can be treated in
the same way.

\ms {\bf 6)} The cup product with $c_r(N^*)$ in the situation
above is always injective. As a consequence the Gysin sequence
above for the pair $(\ov{\cL}_{\le p};\ov{\cL}_p)$ splits
completely into short exact sequences
$$
0\to H^{\bullet -2c}(\ov{\cL}_p) \to H^\bullet(\ov{\cL}_{\le
p})\to H^\bullet(\ov{\cL}_{<p})\to 0.
$$
The $\bbbz/l^n\bbbz$-modules $H^\bullet(\ov{\cL}_{\le
p};\bbbz/l^n\bbbz)$ are free of finite rank.

\ms {\bf Proof.} Using induction over $n$ and the order of the
Harder-Narasimhan polygon, it suffices to treat the case $n=1$.
Using 3), one is reduced to the corresponding statement for
$\ov{\cL}(X;n;d)_{ss}$ and  products of these resp. to the
``$0$''-versions of these. Using the usual geometric invariant
theory description of the moduli space of vector bundles, such a
product can be written as a quotient stack $[G\setminus
\prod\limits^l_{i=1}Gr_m(n_i;d_i)]$ of a product of Grassmannians
by the action of a connected algebraic group $G$. By (\cite{Be},
Theorem 1.4.3), the spectral sequence describing this quotient
degenerates and in particular the cohomology of the stack injects
into the cohomology of the corresponding quotient stack
$[T\setminus \prod\limits^l_{i=1}Gr_m(n_i;d_i)]$, on which it is
easy to check the injectivity of the cup product above. A
different proof can be derived also from \cite{Bi-Ghi-Lat}, in
particular using Prop. 4.2, Prop. 7.1.\qed

\ms Using a Mittag-Leffler type argument one obtains now
$$H^\bullet(\ov{\cL};\bbbq_l)=\lim\limits_{\overleftarrow{(p)}}H^\bullet(\ov{\cL}_{\le
p};\bbbq_l)$$ and it follows

{\proposition  The Gysin sequence for the pair $(\ov{\cL}_{\le p},
\ov{\cL}_p)$ splits into the short exact sequences
$$ 0 \to H^{\ast-2r} (\ov{\cL}_p; \ov{\bbbq}_l
(r))\mathop{\to}\limits^{i_{\ast}}  H^{\ast} (\ov{\cL}_{\le p};
\ov{\bbbq}_l) \to H^{\ast} (\ov{\cL}_{< p}; \ov{\bbbq}_l ) \to 0,
 $$
where $r$ is the codimension of the stack $\ov{\cL}_p$ in
$\ov{\cL}_{\le p}$.}

\ms {\bf 7)} Now we have the following result concerning the
Poincar\'e series of the cohomology ring of the moduli stack
$\ov{\cL}$

\ms {\proposition The Poincar\'e series of the cohomology ring
$H^{\ast}(\ov{\cL}; \ov{\bbbq}_l)$ is given as}
$$
P (\ov{\cL}; t) = \displaystyle\frac{\prod\limits^n_{i=1} (
1+t^{2i-1})^{2g}}{\prod\limits^n_{i=2}
(1+t^{2i})\prod\limits^n_{i=2}(1-t^{2i-2})}.
$$

\ms {\bf Proof.} This follows now directly from Proposition 2.1
using \cite{Bi-Ghi-Lat}, Prop. 4.2, Prop. 8.1, and Prop. 10.1.. An
alternative, more arithmetic proof is given in \cite{Hei}, Satz
2.2.6, but with the slight modification, that we are considering
the case with trivial determinant. This proof uses again
Proposition 2.1., the Lefschetz trace formula (see \cite{Be}) and
the Tamagawa number computations from \cite{Ha-Na}. \qed

\ms {\bf 8)} The cohomology of $\ov{\cL}$ can be obtained now as
follows: we consider on  $\ov{X} \times \ov{\cL}$ the universal
bundle $E^{\univ}$ and its Chern classes
$$
c_i (E^{\univ}) \in H^{2i} (\ov{X} \times \ov{\cL}; \ov{\bbbq}_l)
\quad (i \ge 2).
$$

\ms {\bf Remark.} The first Chern class $c_1 (E^{\univ})$ vanishes
because we have by definition the trivialisation $\delta^{\univ}:
\det(E^{\univ}) \stackrel{\sim}{\to} \cO_{\ov{X}\times \ov{\cL}}$.

\ms We apply now K\"unneth's theorem, which holds also in the
situation here, because we can view $\ov{\cL}$ as a simplicial
scheme using a representation of $\ov{\cL}$ as mentioned above and
therefore also the product $\ov{X} \times \ov{\cL}$ as a
simplicial scheme. Then we can apply K\"unneth's theorem
degreewise for the simplicial scheme representing $\ov{X} \times
\ov{\cL}$ to obtain it for $\ov{X} \times \ov{\cL}$ by standard
simplicial techniques. Therefore we can decompose the Chern
classes above as follows:
$$c_i(E^{\univ}) = 1 \otimes c_i + \sum\limits^{2g}_{j=1} \alpha_j \otimes
a^{(j)}_i + [X] \otimes b_{i-1}.$$

Here we have, for $i \ge 2$, $c_i \in H^{2i}(\ov{\cL};
\ov{\bbbq}_l)$, $a^{(j)}_i \in H^{2i-1} (\ov{\cL}; \ov{\bbbq}_l)$
for $j=1, \ldots, 2g$ and $b_{i-1} \in H^{2(i-1)} (\ov{\cL};
\ov{\bbbq}_l)$.

\ms {\proposition There is an isomorphism of graded
$\ov{\bbbq}_l$-algebras
$$
H^{\ast} (\ov{\cL}; \ov{\bbbq}_l) \cong \ov{\bbbq}_l [c_2, \ldots,
c_n; b_1, \ldots, b_{n-1}] \otimes \bigotimes\limits^n_{i=1}
\Lambda [a^{(1)}_i, \ldots, a^{(2g)}_i],
$$
where $\bbbq_l[c_2,\dots,c_n; b_1,\ldots, b_{n-1}]$ is the
(graded) polynomial algebra in the generaters $c_2, \ldots, c_n,
b_1 \ldots b_{n-1}$ and $\Lambda [a^{(1)}, \ldots, a^{(2g)}_i]$ is
an exterior algebra with generators $a^{(1)}_i, \ldots,
a^{(2g)}_i$ in degree $(2i-1)$.}

\ms {\bf Proof.} {The proof can be found in \cite{Hei}, Satz 2.2.8
or \cite{Bi-Ghi-Lat}, again of course with the slight
modification, that we are considering the case of trivial
determinant. The strategy of the proof can be outlined as follows:
first, upon restriction to the closed substack of $\ov{\cL}$,
consisting of vector bundles, which are direct sums of line
bundles, it follows that the canonical map
$$
\ov{\bbbq}_l [c_1, \ldots, c_n ; b_1,\ldots, b_{n-1}] \otimes
\bigotimes\limits^n_{i=1} \Lambda [a^{(1)}_i, \ldots, a^{(2g)}_i]
\to H^{\ast}(\ov{\cL}; \ov{\bbbq}_l),
$$
defined by sending the generators $c_i, b_i, a^{(j)}_i$ to the
corresponding cohomology classes in $H^{\ast} (\ov{\cL};
\ov{\bbbq}_l)$, is injective. The surjectivity follows then from
comparing the Poincar\'e polynomials of the two graded algebras
using Proposition 2.2.} \qed

\ms We will study now the action induced by the morphism
$\varphi:\bar{\cL}\to\bar{\cL}$ of stacks on cohomology. By
Proposition 1.6. we have the following isomorphism of vector
bundles on $\bar{X}\times \bar{\cL}$,
$$
(\ov{F}_X\times \id_{\bar{\cL}})(E^{{\rm univ}})\cong
(\id_{\bar{X}}\times \varphi)^*(E^{{\rm univ}})
$$
This induces the equalities of Chern classes $(i \ge 2$)
$$
c_i((\ov{F}_X\times \id_{\bar{\cL}})^*(E^{{\rm
univ}}))=c_i((\id_{\bar{X}}\times \varphi)^*(E^{{\rm univ}}))
$$
Using the functoriality of Chern classes, we obtain therefore
$$
(F_{\bar{X}}\times\id_{\bar{\cL}})^*(c_i(E^{{\rm
univ}}))=(\id_{\bar{X}}\times\varphi)^*(c_i(E^{{\rm univ}}))
$$
From above we have the formula
$$
c_i(E^{{\rm univ}})=1\otimes c_i+\sum^{2g}_{j=1}\alpha_j\otimes
a^{(j)}_i+[X]\otimes b_{i-1}
$$
for $i \ge 2$. Therefore we obtain the equality \beqn & & 1\otimes
c_i+\sum\limits^{2g}_{j=1}\lambda_j\alpha_j\otimes
a^{(j)}_i+q[X]\otimes b_{i-1}\\
 &=&1\otimes \varphi^*(c_i)+\sum\limits^{2g}_{j=1}\alpha_j\otimes
\varphi^*(a^{(j)}_i)+[X]\otimes \varphi^*(b_{i-1}) \eeqn

This implies the following equalities
$$
\ba{ll}
\varphi^*(c_i)=c_i& (i\ge 1)\\
\varphi^*(a^{(j)}_i)=\lambda_ja^{(j)}_i  &(i\ge 1,\, j=1,\ldots,2g)\\
\varphi^*(b_i)=q\, b_i&(i\ge 1) \ea
$$
and we have proved: \ms {\proposition The morphism
$\varphi:\bar{\cL}\to\bar{\cL}$ of the stack of vector bundles of
rank $n$ on $\bar{X}$, given by the pullback operation induced by
$\ov{F}_X:\bar{X}\to \bar{X}$, acts on the generating cohomology
classes $c_i,\, a^{(j)}_i,\, b_i$ of $H^*(\bar{\cL};\ov{\bbbq}_l)$
as:}
$$
\ba{ll}
\varphi^*(c_i)=c_i& (i\ge 1)\\
\varphi^*(a^{(j)}_i)=\lambda_j\, a^{(j)}_i  &(i\ge 1,\,
j=1,\ldots,2g)\\\varphi^*(b_i)=q\, b_i&(i\ge 1) \ea
$$

\ms {\bf Remark.} As mentioned earlier, it is somewhat surprising
that one can determine the action of $\varphi$ on the cohomology
of $\cL$ in such an explicit way, because otherwise the nature of
$\varphi$ is mysterious.

\ms We will use now similar considerations to study the action of
the genuine Frobenius endomorphism $\ov{F}_{\cL}=F_{\cL}\times
\id_{\ov{\bbbf}_q/\bbbf_q}$
on the cohomology $H^*(\bar{\cL}; \ov{\bbbq}_l)$.\\
To proceed further, we consider the classifying stacks of all rank
$n$ vector bundles $BGL(n)/\bbbf_q$ and
$BGL(n)/\bbbf_q\times\ov{\bbbf}_q$. As is well known, (see
\cite{Be}, Th. 2.3.2.) one has
$$
H^*(BGL(n)\times\ov{\bbbf}_q; \ov{\bbbq}_l)\simeq
\ov{\bbbq}_l[c_1,\ldots,c_n]
$$
and the geometric Frobenius $\ov{F}_{BGL(n)/\bbbf_q}$ acts as
$$
(\ov{F}_{BGL(n)/\bbbf_q})^*(c_i)=q^ic_i\, \, \, \mbox{\rm for}\,
\, i\ge 1.
$$
Using this we can obtain more information about the Frobenius
action on the cohomology $H^*(\bar{\cL}; \ov{\bbbq}_l)$. We have
$$
(\ov{F}_X\times
\id_{\bar{\cL}})^*(\id_{\bar{X}}\times\ov{F}_{\cL})^*(c_i(E^{{\rm
univ}}))\cong(\ov{F}_{X\times \cL})^*(c_i(E^{{\rm univ}}))\cong
(\ov{F}_{X\times \cL})^*c_i(u^*\tilde{E}^{{\rm univ}}),
$$
where $\tilde{E}^{{\rm univ}}$ denotes the universal vector bundle
on the classifying stack $BGL(n)/\bbbf_q$ and $u:X\times\cL\to
BGL(n)$ denotes the classifying morphism for $E^{{\rm univ}}$ on
$(X\times \cL)$, such that
$$
u^*(\tilde{E}^{\univ})\cong E^{\univ}
$$
with a canonical isomorphism. Furthermore we have the commutative
diagram
$$
\ba{ccc}
X\times \cL & \vrpf{10}{u}& BGL(n)/\bbbf_q\qquad\\
\vupflr{10}{F_{X\times\cL}}{} & & \vupflr{10}{}{F_{BGL(n)/\bbbf_q}} \\
X\times\cL & \vrpf{10}{u} & BGL(n)/\bbbf_q \ea
$$
and similarly after extension from $\bbbf_q$ to $\ov{\bbbf}_q$.
Therefore we conclude \beqn
(\ov{F}_{X\times\cL})^*c_i(\bar{u}^*(\tilde{E}^{{\rm
univ}}))&=&(\ov{F}_{X\times \cL})^*(\bar{u})^*c_i(\tilde{E}^{{\rm univ}})\\
 &=& (\bar{u})^*(\ov{F}_{BGL(n)/\bbbf_q})^*(c_i(\tilde{E}^{{\rm univ}}))\\
 &=& q^i(\bar{u})^*(c_i(\tilde{E}^{{\rm univ}}))\\
&=& q^i(c_i((\bar{u})^*\tilde{E}^{{\rm univ}}))\\
&=& q^ic_i(E^{{\rm univ}}). \eeqn On the other hand we obtain also
\beqn (\ov{F}_X\times \id_{\bar{\cL}})^*(\id_{\bar{X}}\times
\ov{F}_{\bar{\cL}})^*(c_i(E^{{\rm univ}})) =(\id_{\bar{X}}\times
F_{\bar{\cL}})^*(\ov{F}_X\times\id_{\bar{\cL}})^*(c_i(E^{{\rm univ}}))\\
=(\id_{\bar{X}}\times
F_{\bar{\cL}})^*(\ov{F}_X\times\id_{\bar{\cL}})^*(1\otimes
c_i+\sum^{2g}_{j=1}\alpha_j\otimes a^{(j)}_i+[X]\otimes b_{i-1})\\
=(\id_{\bar{X}}\times F_{\bar{\cL}})^*(1\otimes
c_i+\sum^{2g}_{j=1}\lambda_j\alpha_j\otimes a^{(j)}_i+q[X]\otimes
b_{i-1}) \eeqn This implies immediately the following fundamental
proposition

\ms {\proposition The action of the geometric Frobenius
$\ov{F}_{\cL}$ on the generating classes of the cohomology
$H^*(\bar{\cL}; \ov{\bbbq}_l)$ of the stack $\cL$ of rank $n$
vector bundles on $X$ is given as} \beqn
\ov{F}^*_\cL(c_i)&=&q^ic_i,\,\, \,  (i\ge 2)\\
\ov{F}^*_\cL(b_i)&=&q^{i-1}b_i,\, \, (i\ge 1)\\
\ov{F}^*_\cL(a^{(j)}_i)&=&\lambda^{-1}_jq^ia^{(j)}_i\, \, (i\ge
2,\, j=1,\ldots 2g). \eeqn

\ms {\corollary The action of the geometric Frobenius
$F_{\bar{\cL}}$ on the cohomology of the open substack
$\bar{\cL}_{{\rm semistable}}$ of semistable bundles is
semisimple.}

\ms {\bf Proof.} {The restriction morphism $H^*(\bar{\cL};
\ov{\bbbq}_l)\to H^*(\bar{\cL}_{{\rm semistable}}; \ov{\bbbq}_l)$
is surjective. As the cohomology $H^*(\bar{\cL}; \ov{\bbbq}_l)$ is
a semisimple $\ov{\bbbq}_l[F_{\bar{\cL}}]$-module, the same is
also true for $H^*(\bar{\cL}_{{\rm semistable}}; \ov{\bbbq}_l).
$\qed}

\ms {\bf Remark.} The semisimplicity of the Frobenius here follows
basically from the semisimplicity of the action of the Frobenius
on the cohomology of a curve, using \cite{Mu}, p. 203. There the
result is proven for abelian varieties over finite fields. As the
first cohomology of a curve can be identified with that of its
Jacobian, the desired result for curves follows immediately.

\ms We have considered above only the case of bundles with trivial
determinant. The same technique can be applied however for the
stack of bundles of rank $n$ and degree $d$, where $\gcd (n,d)=1$,
i.e. $n$ and $d$ are coprime numbers.

\ms {\proposition Suppose $\gcd (n,d)=1$ as above. Then the action
of the Galois group $\Gal(\ov{\bbbf}_q/\bbbf_q)$ resp. the
arithmetic Frobenius $\psi$ on the cohomology
$H^*(M(n,d);\ov{\bbbq}_l)$ of the coarse moduli space $M(n,d)$ of
stable vector bundles of rank $n$ and degree $d$ is semisimple.}

\ms {\bf Proof.} {Because $\gcd (n,d)=1$, we have
$\tilde{\cL}(n;d)_{{\rm semistable}}=\tilde{\cL}(n,d)_{{\rm
stable}}$. Furthermore there is a canonical morphism
$$
\tilde{\cL}(n;d)_{{\rm stable}}\to M(n,d),
$$
given functorially for any scheme $S/\bbbf_q$ by the obvious
functor
$$
\cL(n,d)_{{\rm stable}}(S)\to M(n,d)(S),
$$
where $M(n,d)(S)$ denotes the discrete category of isomorphism
classes of vector bundles of rank $n$ and degree $d$ (pointwise in
$S$) on
$X\times_{\bbbf_q} S$.\\
The Leray spectral sequence for this morphism is degenerating,
because we have an equality of Poincar\'e polynomials
$$
P(\bar{\cL}(n;d)_{{\rm stable}};t)=P(\ov{M}(n,d);t)\cdot
P(B\bbbg_m;t),
$$
where the morphism $\bar{\cL}_{{\rm stable}}\to \ov{M}(n;d)$ is a
fibration
with fibre the classifying stack $B\bbbg_m$.\\
Therefore $H^*(\ov{M}(n;d); \ov{\bbbq}_l)$ is a submodule of
$H^*(\bar{\cL}_{{\rm stable}}; \ov{\bbbq}_l)$. As the last one is
a semisimple module over $\ov{\bbbq}_l[\ov{F}_{\cL}]$,
$H^*(\ov{M}(n;d); \ov{\bbbq}_l)$ is a semisimple module over
$\ov{\bbbq}_l[\ov{F}_{M(n;d)}]$, or equivalently the Frobenius
$\psi$ is acting as a semisimple endomorphism.}\qed

\ms {\proposition Identifying the algebraic closure $\ov{\bbbq}\,
_l$ with the field $\bbbc$ of complex numbers, the expression for
the formal trace
$$
\tr(\varphi^r\times \psi^s;\, H^*(\bar{\cL};
\ov{\bbbq}_l))=\sum_{i\ge 0}(-1)^i\tr(\varphi^r\times
\psi^s;H^i(\bar{\cL};\ov{\bbbq}_l))
$$
is absolutely convergent for $s > r$. }

\ms {\bf Proof.} {This is just an exercise in summing up geometric
series. The formal series to be considered is given as
$$
(\prod\limits^n_{i=2}(\sum^\infty_{n=0}q^{-(si)}))\cdot(\prod\limits^{n-1}_{k=1}\,
\,
\sum^\infty_{m=0}q^{m(r-ks)})\cdot(\prod\limits^{2g}_{j=1}\prod\limits^{n}_{i=1}(1+|\lambda_j|^{r+s}q^{-is}))
$$
which implies the statement.} \qed

\ms {\bf Remark.} The question is now of course, if is there a
kind of Lefschetz trace formula in this general situation.

\ms This is true at least in the case $r=0$, $s=1$, which was
considered in \cite{Be}, \cite{Hei}. The following trace formula
holds
$$
q^{\dim(\cL)}\sum_{i\ge 0}(-1)^i\tr(\psi; H^i(\bar{\cL};
\ov{\bbbq}_l))=\sum_{[E]\in \cL(k)}\frac{1}{|\Aut^{(0)}(E)|},
$$
where
$$
\Aut^{(0)}(E):=\{\alpha\in \Aut(E)(k)\mid \det(\alpha)=1\}
$$
and $\dim(\cL)=n^2(g-1)+1$.

\ms But as mentioned already in the introduction, the situation in
the cases $r=0$ resp.\ $ r > 0$ is rather different, because the
proof in the case $r=0$ works by studying $\bar{\cL}$ and the
traces above using the filtration by the $\bar{\cL}_{\le p}$ as
done earlier. On the other hand the morphism
$\varphi:\bar{\cL}\to\bar{\cL}$ does not respect these open
substacks.

\section{Some complements and examples}

\ms {\bf Part 1: Vector bundles of rank 2 on the projective line}

\ms It seems worthwhile to study in detail the easiest example,
namely the case of the projective line $X=\bbbp^1$ and vector
bundles of rank $n=2$ with trivial determinant. We specialize our
computations from section 2 to this case. \ms By Proposition 2.3,
we have in this special situation
$$
H^*(\bar{\cL}(\bbbp^1;2); \ov{\bbbq}_l)\cong
\ov{\bbbq}_l[c_2,b_1],
$$
where
$$
b_1\in H^2(\bar{\cL}; \ov{\bbbq}_l),\, c_2\in H^4(\bar{\cL};
\ov{\bbbq}_l)
$$

For the formal trace we obtain \beqn \tr(\varphi^r\times \psi^s;\,
H^*(\bar{\cL}(\bbbp^1;2); \ov{\bbbq}_l) &=& \sum_{i\ge
0}(-1)^i\tr(\varphi^r\times \psi^s;H^i(\bar{\cL}(\bbbp^1;2);
\ov{\bbbq}_l)\\
 &=&(\sum^\infty_{m=0}q^{-2sm})(\sum^\infty_{m=0}q^{(r-s)m}),
\eeqn which is convergent for $r < s$. \ms We obtain therefore for
$r < s$
$$
\tr(\varphi^r\times \psi^s; H^*(\bar{\cL}(\bbbp^1;2);
\ov{\bbbq}_l))=(1-q^{-2s})^{-1}(1-q^{r-s})^{-1}
$$

\ms {\bf Remark.} In this special case the fixed points under the
action of $(\varphi^r\times\psi^s)$ can be computed directly. The
naive fixed point set can be computed as
$$
(\bar{\cL})^{\varphi^r\times\psi^s}=(\bar{\cL}_{\rm
semistable})^{\varphi^r\times
\psi^s}=(BSL(2)/\bbbf_q\times\ov{\bbbf}_q)^{\varphi^r\times\psi^s}
$$
Here we have used the isomorphism of stacks
$$
BSL(2)/\bbbf_q\to \cL(\bbbp^1;2)_{\rm semistable},
$$
given functorially for the categories of $S$-valued points,
$S/\bbbf_q$ a scheme, by
$$
E\longmapsto \pr^*_2(E),
$$
where $E$ is a vector bundle of rank two over $S$ with trivial
determinant and $\pr^*_2(E)$ is the pullback to
$\bbbp^1\times_{\bbbf_q}S$ via the projection
$$
\pr_2:\bbbp^1\times_{\bbbf_q}S \to S
$$
This implies immediately, that $\varphi$ acts as the identity on
$\cL_{\rm semistable}$. But then we obtain
$$
(\bar{\cL})^{\varphi^r\times\psi^s}=(BSL(2)/\bbbf_q\times\ov{\bbbf}_q)^{\psi^s}
$$
This does not depend on $r$ at all, which obviously does not fit
with the formula obtained above.

\bigskip
{\bf Part 2: Affine Grassmannians and gauge groups}

\ms In this second part of section 3 we give a description of the
moduli stack of vector bundles as a double quotient similar to the
well known adelic description of vector bundles on a curve. These
considerations are useful for our goal to make sense of a general
trace formula, in order to do our computations of section 2 for an
appropriate compactification of the stack $\bar{\cL}$ or $\cL$.

\ms For the following we consider a smooth projective curve $X$
over $\bbbf_q$ with a closed $\bbbf_q$-rational point $\infty\in
X$ (whose existence we assume for simplicity). $\cO_{X,\infty}$ is
the local ring of $X$ at $\infty$, $\cO_\infty$ its completion,
$F_\infty$ the completion of the function field $F=\bbbf_q(X)$ at
$\infty$. $X^{(\infty)}=X\setminus\{\infty\}$ is an open subscheme
of $X$.

\ms Let $\tilde{\cL}(X;n)$ denote the algebraic stack of vector
bundles of rank $n$ on $X$, similarly
$\tilde{\cL}(X^{(\infty)};n)$ denotes the stack of vector bundles
of rank $n$ on $X^{(\infty)}$. Finally, let $\cL(X;n)$ resp.\
$\cL(X^{(\infty)};n)$ denote the corresponding stacks with
trivialized determinant.

\ms There is an obvious morphism $r$ of stacks
$$
r:\cL(X;n)\to \cL(X^{(\infty)};n)
$$
given by the restriction functor
$$
E\in \cL(X;n)(S)\longmapsto (E\mid
X^{(\infty)}\times_{\bbbf_q}S)\in\cL(X^{(\infty)};n)(S),
$$
where $S/\bbbf_q$ is any scheme over $\bbbf_q$.

\ms The fibre over a closed point $(E^{(0)};\delta^0)\in
\cL(X^{(\infty)};n)$ can be described as the stack, given
functorially by $S/\bbbf_q\mapsto \underline{\underline{\rm
Gr}}(E^{(0)})(S)$, where  $\underline{\underline{\rm
Gr}}(E^{(0)})(S)$ denotes the category of triples
$((E,\delta),n)$, where $(E,\delta)\in \cL(X;n)(S)$ that is, $E$
is a rank $n$ vector bundle over $X\times_{\bbbf_q}S$, $\delta
:\det(E)\stackrel{\sim}{\to}\cO_{X\times_{\bbbf_q}S}$ a
trivialisation of the associated determinant bundle and $u$ is an
isomorphism of vector bundles
$$
u:E\mid_{X^{(\infty)}\times_{\bbbf_q}S}\widetilde{\longrightarrow}
\pr^*_2(E^{(0)}),
$$
where $\pr_2:X^{(\infty)}\times_{\bbbf_q}S\to S$ is the projection
on the second factor. Additionally, the induced determinantal
isomorphism $\det(u)$ has to make the obvious diagram of
determinant bundles and trivialisations commutative.

\ms {\bf Remark.} For simplicity we have written
$\underline{\underline{\rm Gr}}(E^{(0)})$ instead of
$\underline{\underline{\rm Gr}}((E^{(0)};\delta^{(0)}))$.

\ms As is known (see for example \cite{Gai}), there is a purely
local description of $\underline{\underline{\rm Gr}}(E^{(0)})$ as
follows:

\ms Consider the group-valued functors \beqn \underline{\mbox{\rm
affine
schemes}}/\bbbf_q&\longrightarrow&\underline{\mbox{\rm groups}}\\
S&\longmapsto&\SL(n;\cO_S[[t_\infty]]),\\
\mbox{\rm resp.}\qquad S&\longmapsto& \SL(n;\cO_S((t_\infty)))
\eeqn

Then the first functor is represented by an algebraic group (not
of finite type) and the second functor is represented by an
ind-algebraic group. We denote these two objects somewhat
informally as $\underline{\underline{{\rm SL}}}(n;\cO_\infty)$
resp. $\underline{\underline{{\rm SL}}}(n;F_\infty)$.

\ms The quotient
$$
\underline{\underline{\SL}}(n;F_\infty)/\underline{\underline{\SL}}(n;\cO_\infty)
$$
exists as a stack (or even an ind-algebraic scheme) and is
isomorphic to the stack $\underline{\underline{\Gr}}(E^{(0)})$.
(see \cite{Gai} and \cite{Fal2}.)

\bigskip
Choosing a pair $(E^{(0)},\delta^{(0)})\in
\cL(X^{(\infty)};n)(\bbbf_q)$, for example the trivial bundle, we
have the following diagram of stacks
$$
\underline{\underline{\Gr}}(E^{(0)})\to \cL(X;n)\to
\cL(X^{(\infty)};n),
$$
where $\underline{\underline{\Gr}}(E^{(0)})\to \cL(X;n)$ is the
morphism of stacks, given functorially by
$$
((E,\delta),u)\mapsto (E,\delta)
$$
using the notations from above.

\ms Now we indicate the computation of the cohomologies in this
case.

\ms {\proposition The filtration of the stack $\cL(X;n)$ by the
open substacks $\cL(X;n)_{\le p}$ induces a corresponding
filtration $\underline{\underline{\Gr}}(E^{(0)})_{\le p}$ of
$\underline{\underline{\Gr}}(E^{(0)})$.}

\ms This can be used to compute the cohomology of the infinite
Grassmannian $\underline{\underline{\Gr}}(E^{(0)})$ as follows:

\ms {\proposition The cohomology ring
$H^*(\underline{\underline{\Gr}}(E^{(0)})\times \ov{\bbbf}_q;
\ov{\bbbq}_l)$ is given as $\ov{\bbbq_l}[b_1,\ldots,b_{n-1}]$
where the $b_i$ $(i=1,\ldots,n-1)$ are the restrictions of the
corresponding classes $b_i$ in the cohomology ring $H^*(\bar{\cL};
\ov{\bbbq}_l)$.}

\ms {\bf Proof.} {We can consider the universal object
$$(E^{{\rm univ}},\delta^{{\rm univ}}),u:E^{{\rm
univ}}\mid_{X^{(\infty)}\times_{\bbbf_q}\underline{\underline{\Gr}}(E^{(0)})}\widetilde{\longrightarrow}
\pr^*_2(\cO^n_{X^{(\infty)}})
$$
over $X\times_{\bbbf_q}\underline{\underline{\Gr}}(E^{(0)})$. We
can write down (for example by pullback) the Chern classes of
$E^{{\rm univ}}$ as
$$
c_i(E^{{\rm univ}})=1\otimes c_i+\sum^{2g}_{j=1}\alpha_j\otimes
a^{(j)}_i+[X]\otimes b_{i-1}
$$
for $i\ge 1$ in
$H^{2i}((X\times_{\bbbf_q}\underline{\underline{\Gr}}(E^{(0)}))\times_{\bbbf_q}\ov{\bbbf}_q;
\ov{\bbbq}_l)$.

Here the $c_i,b_i,a^{(j)}_i$ are just the restrictions of the
corresponding classes under the morphism
$\underline{\underline{\Gr}}(E^{(0)})\to\cL$ considered above. On
the other hand we have the identity
$$
c_i(E^{{\rm univ}})\mid (X^{(\infty)}\times
\underline{\underline{\Gr}}(E^0))\times
\ov{\bbbf}_q=\varphi^*(\pr^*_1(c_i(E^{(0)}))
$$
by the very definition of $(E^{{\rm univ}},\varphi^{{\rm univ}})$.\\
But $c_i(E^{(0)})=0$ for $i\ge 1$, therefore we obtain for the
restrictions \beqn
c_i\mid \underline{\underline{\Gr}}(E^0)&=&0\qquad (i\ge 1)\\
a^{(j)}_i\mid \underline{\underline{\Gr}}(E^0)&=&0\qquad (i\ge
1,\, j=1,\ldots,2g) \eeqn It remains to show that
$b_1,\ldots,b_{n-1}$ are generating elements of the cohomology and
are polynomially independent. This can be done using the
stratification mentioned above in a similar treatment as in the
corresponding statements for the cohomology of the stack
$\bar{\cL}$.\qed}

\ms {\bf Remark.} In a similar way, using a yet to be defined
version of $l$-adic cohomology appropriate for the stack
$\cL(X^{(\infty)}\times\ov{\bbbf}_q;n)$, we could consider the
morphism $r:\cL(X;n)\to \cL(X^{(\infty)};n)$, given as restriction
of bundles to $X^{(\infty)}$. For the universal bundles we have
the pullback situation
$$
r^*(E^{{\rm univ}}({\rm over}\,
X^{(\infty)}\times\cL(X^{(\infty)};n))=E^{{\rm univ}}({\rm over}\,
X\times \cL(X;n))
$$
Again for $E^{{\rm univ}}$ over $X^{(\infty)}\times
\cL(X^{(\infty)};n)$, we could consider the Chern classes and then
obtain the following description
$$
c_i(E^{{\rm univ}})=1\otimes c_i+\sum^{2g}_{j=1}\alpha_j\otimes
a^{(i)}_j,\, (1\le i\le n)
$$
because here we have $[X]|H^2(\bar{X}_\infty)=0$.

\ms Then the cohomology
$H^*(\cL(X^{(\infty)};n)\times\ov{\bbbf}_q;\ov{\bbbq}_l)$ is given
as
$$
\ov{\bbbq}_l[c_1,\ldots,c_n]\otimes
\bigotimes\limits^n_{i=1}\Lambda[a^{(1)}_i,\ldots,a^{(2g)}_i]
$$
in the sense of Proposition 2.3 and the Leray spectral sequence
for the morphism of stacks
$$
\bar{r}:\bar{\cL}(X;n)\to \bar{\cL}(X^{(\infty)};n)
$$
would degenerate.

\bigskip
We would get this description along the following lines. As the
cohomology classes
$$
\{c_i(\bar{\cL}(X)), a^{(j)}_i(\bar{\cL}(X))\mid i=1,\ldots,n,\,
j=1,\ldots,2g\}
$$
$$
\{b_i(\bar{\cL}(X))\mid i=1,\ldots,n-1\}
$$
are independent in the sense of Proposition 2.3, the same is also
true for their preimages
$$
\{c_i(\bar{\cL}(X^{(\infty)}));a^{(j)}_i(\bar{\cL}(X^{(\infty)}))\mid
i=1,\ldots,n,\,{\rm and}\, \,  j=1,\ldots,2g\}.
$$
On the other hand comparing the Poincar\'e polynomials we would
get
$$
P(\cL(X);t)=P(\underline{\underline{\Gr}}(E^{(0)};t)P(\cL(X^{(\infty)});t)
$$
From this and the corresponding computations of
$H^*(\bar{\cL}(X);\ov{\bbbq}_l)$ and the cohomology
$H^*(\underline{\underline{\Gr}}(E^{(0)});\ov{\bbbq}_l)$ we would
then obtain the desired description of the cohomology rings and
the Leray spectral sequence.

\end{document}